\documentclass[12pt,a4paper]{article}
\usepackage{geometry}
\usepackage{graphicx}
\usepackage{hyperref}
\usepackage{url}
\usepackage{tikz}
\usetikzlibrary{calc}
\usepackage{asymptote}
\usepackage{amsmath,amssymb,amsthm}

\theoremstyle{plain}

\theoremstyle{definition}

\title{Generalization of Marion's Theorem:\\ Volumes of Central Polytopes Obtained by Trisecting the Edges of Simplices}

\usepackage{authblk}
\author{Yu.~V.~Kazakov\thanks{E-mail: \texttt{yu.v.kazakov@gmail.com}}}
\affil{%
    Reshetnev Siberian State University of Science and Technology(retired) \\
    \small Address: 31, Krasnoyarskii Rabochii Ave., 82, Mira Ave., 660037, Krasnoyarsk, Russia
}
\date{}

\begin{document}
\maketitle

\begin{abstract}
The classical Marion's theorem states that the area of the central hexagon obtained by dividing each side of a triangle into three equal parts and connecting the interior division points is exactly $1/10$ of the area of the original triangle.
In this paper, the construction is extended to an arbitrary $n$-dimensional simplex: each edge is divided into three equal parts, and through each $n-1$ vertices and the two interior division points of the opposite edge, two hyperplanes are drawn.
It is shown that the resulting inner polytope has a surprisingly simple volume formula: the volume is $$\frac{1}{\binom{2n+1}{n}}$$ times the volume of the original simplex.
The proof uses a probabilistic interpretation with exponential random variables and reduces the geometric problem to a combinatorial binomial sum.
Additionally, the combinatorics of the polytope are discussed: the numbers of vertices, edges, and faces are expressed in closed form, and the complete $f$-vector is described.
A link to an interactive three-dimensional illustration for a regular tetrahedron in GeoGebra \cite{GeoGebra3D} is provided.
For a clearer understanding of the solution approach, Appendix A presents a detailed analysis of the probabilistic method using the triangle as an example, Appendix B presents the technique of applying barycentric coordinates to solve the problem in the three-dimensional case, and Appendix C contains Asymptote \cite{Asymptote} code for visualizing the central deltoidal dodecahedron inside a regular tetrahedron.
\end{abstract}

\section{Introduction}
\begin{figure}[htbp]
\centering
\includegraphics[width=0.45\textwidth]{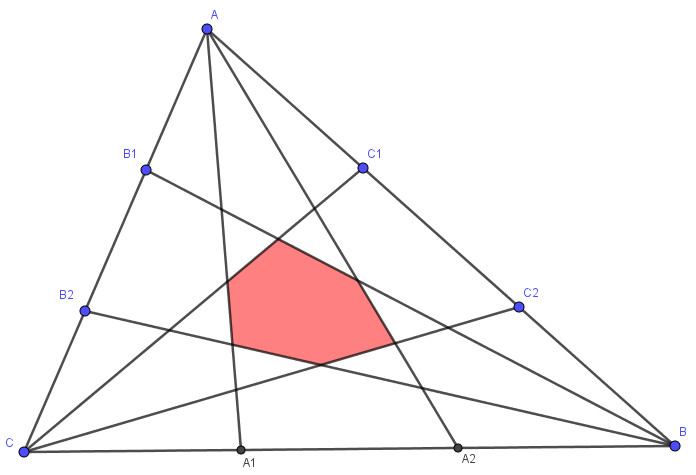}
\caption{Marion's central hexagon in a triangle. Image taken from \cite{DriveImage}.}
\label{fig:marion}
\end{figure}

Marion's theorem~\cite{Marion1993,MathWorldMarion} is an elegant elementary result in plane geometry: in any triangle, the hexagon formed by the interior points dividing each side in the ratio $1:2$ has an area exactly one tenth of the area of the whole triangle (Fig.~\ref{fig:marion}).
Through each vertex of the triangle, two lines are drawn to the two trisection points of the opposite side; the region between these lines in the center of the triangle forms a hexagon. The area fraction $1/10$ of this hexagon is independent of the shape of the triangle.

Generalizations of this construction to higher dimensions are absent in the literature — the author first addressed this problem in 2009 and obtained a result for the case of trisecting the edges of a tetrahedron in $3D$ using Mathcad. In three-dimensional space, an analogous procedure applied to a tetrahedron (each edge is divided into three equal parts, planes are drawn through two vertices and an interior point on the opposite edge — 12 planes in total) yields a central polytope whose volume is $1/35$ of the volume of the tetrahedron~\cite{3dgen}.
A live model demonstrating this construction for a regular tetrahedron is available on the GeoGebra website~\cite{GeoGebra3D}; the central polytope in this case is a deltoidal dodecahedron (a name proposed by the author of this note, since the faces of the polytope are identical deltoids, and the polytope is analogous to the rhombic dodecahedron).
In four-dimensional space, the volume fraction turns out to be $1/126$, where for each set of 3 vertices, two hyperplanes are constructed, each containing one of the two points on the opposite edge (the opposite edge is also divided into three equal parts).

The obtained volume fractions of the inner polytopes —
\[
\frac{1}{10} = \frac{1}{\binom{5}{2}},\quad \frac{1}{35} = \frac{1}{\binom{7}{3}},\quad \frac{1}{126} = \frac{1}{\binom{9}{4}}
\]
— unambiguously point to a general formula:
\[
\frac{\operatorname{Vol}(\text{inner polytope})}{\operatorname{Vol}(\text{simplex})} = \frac{1}{\binom{2n+1}{n}}.
\]
The purpose of this paper is to provide a rigorous proof of this formula for arbitrary dimension $n\ge 1$ and to describe the combinatorial structure of the resulting central polytope.

\section{Construction and Inequalities}
Consider an $n$-dimensional simplex $\Delta^n$. In barycentric coordinates $(x_0,x_1,\dots,x_n)$, it is defined by the conditions
\[
x_i \ge 0,\qquad \sum_{i=0}^{n} x_i = 1.
\]
The ratio of the volume of the inner polytope to the volume of the original simplex is an affine invariant; therefore, without loss of generality, we can use this barycentric representation. The volume of the standard simplex in this normalization is $V_0 = 1/n!$, but the specific volume value is irrelevant, since only the volume ratio enters the final result.

Let $\mathbf{e}_0, \mathbf{e}_1, \dots, \mathbf{e}_n$ be the vertices of the simplex. In barycentric coordinates, the vertex $\mathbf{e}_i$ has coordinates $(x_0,\dots,x_n)$ with $x_i = 1$ and $x_j = 0$ for $j \neq i$.

Choose an arbitrary edge of the simplex corresponding to vertices $\mathbf{e}_i$ and $\mathbf{e}_j$ ($i \neq j$).
The two points dividing this edge internally in the ratios $1:2$ and $2:1$ have barycentric coordinates
\[
P_{ij} = \frac{2}{3}\mathbf{e}_i + \frac{1}{3}\mathbf{e}_j,\qquad
Q_{ij} = \frac{1}{3}\mathbf{e}_i + \frac{2}{3}\mathbf{e}_j.
\]

Consider all $n-1$ vertices other than $\mathbf{e}_i$ and $\mathbf{e}_j$. 
A hyperplane passes through these $n-1$ vertices and the point $P_{ij}$; a similar hyperplane passes through the same $n-1$ vertices and $Q_{ij}$. In barycentric coordinates, these hyperplanes correspond to the equations
\[
x_j = 2x_i \qquad\text{and}\qquad x_i = 2x_j.
\]
The half-spaces containing the central region are described by the inequalities
\[
\frac{1}{2}x_j \le x_i \le 2x_j.
\]

Repeating this construction for each edge of the simplex, we obtain for any ordered pair $0\le i,j\le n$, $i\neq j$, the constraint
\begin{equation}\label{eq:ineq}
x_i \le 2x_j .
\end{equation}

The inner polytope $\mathcal{P}_n$ is defined as the set of points satisfying
\[
x_i \ge 0,\quad \sum_{i=0}^n x_i = 1,\quad x_i \le 2x_j\ \ \text{for all}\ i\neq j.
\]
By symmetry, it suffices to require the inequalities $x_i \le 2x_j$ for all $i\neq j$; the reverse inequalities $x_j \le 2x_i$ then hold automatically, as they are obtained by permuting indices.

The polytope $\mathcal{P}_n$ is the intersection of the original simplex with a system of half-spaces bounded by the hyperplanes
\[
x_i = 2x_j \quad (i\neq j).
\]
As will be shown below, it is the convex hull of $2^{n+1}-2$ points whose barycentric coordinates take only two distinct values differing by a factor of two.

\section{Volume Calculation}
To find the volume of $\mathcal{P}_n$, we employ a probabilistic method.

Let $Y_0,Y_1,\dots,Y_n$ be independent exponentially distributed random variables with mean $1$. Then the random vector
\[
\left(\frac{Y_0}{S},\frac{Y_1}{S},\dots,\frac{Y_n}{S}\right),\qquad S=\sum_{k=0}^n Y_k,
\]
is uniformly distributed on the simplex $\Delta^n$~\cite{Rubinstein}.
Consequently,
\[
\frac{\operatorname{Vol}(\mathcal{P}_n)}{\operatorname{Vol}(\Delta^n)}
= \mathbb{P}\bigl( \text{the conditions } x_i \le 2x_j \text{ hold} \bigr)
= \mathbb{P}\bigl( Y_i \le 2Y_j \text{ for all } i,j \bigr).
\]

The event $Y_i \le 2Y_j$ for all $i,j$ means that the ratio of the minimum to the maximum among $Y_k$ is at least $1/2$.
By symmetry, we fix which variable is the smallest.
Let $Y_0$ be the minimum; then the condition becomes $Y_k \le 2Y_0$ for all $k=1,\dots,n$ (since $Y_0\le Y_k$ automatically implies $Y_0 \le 2Y_k$).
The probability that $Y_0$ is the minimum and that $Y_k \le 2Y_0$ for $k\ge 1$ is
\[
\int_0^\infty e^{-y_0} \left( \int_{y_0}^{2y_0} e^{-y}\,dy \right)^{\!n} \,dy_0
= \int_0^\infty e^{-y_0} \bigl(e^{-y_0} - e^{-2y_0}\bigr)^n \,dy_0.
\]
Multiplying by $n+1$ (the number of choices for which variable is the minimum), we obtain
\[
\frac{\operatorname{Vol}(\mathcal{P}_n)}{V_0}
= (n+1) \int_0^\infty e^{-y_0} \bigl( e^{-y_0} - e^{-2y_0} \bigr)^n \,dy_0
= (n+1) \int_0^\infty e^{-(n+1)y_0} \bigl( 1 - e^{-y_0} \bigr)^n \,dy_0.
\]

Expanding $(1 - e^{-y_0})^n$ using the binomial theorem:
\[
(1 - e^{-y_0})^n = \sum_{k=0}^n \binom{n}{k} (-1)^k e^{-k y_0}.
\]
The integral becomes a sum of elementary exponential integrals:
\begin{align*}
\frac{\operatorname{Vol}(\mathcal{P}_n)}{V_0}
&= (n+1) \sum_{k=0}^n \binom{n}{k} (-1)^k \int_0^\infty e^{-(n+1+k)y_0}\,dy_0 \\
&= (n+1) \sum_{k=0}^n \binom{n}{k} \frac{(-1)^k}{n+1+k}.
\end{align*}

A known combinatorial identity (see, e.g., \cite{GKP}) states that
\[
\sum_{k=0}^n \binom{n}{k} \frac{(-1)^k}{m+k} = \frac{1}{m \binom{m+n}{n}} \qquad (m>0).
\]
Setting $m = n+1$, we find
\[
\sum_{k=0}^n \binom{n}{k} \frac{(-1)^k}{n+1+k} = \frac{1}{(n+1)\binom{2n+1}{n}}.
\]
Multiplying by $n+1$, we arrive at the main result
\begin{equation}\label{eq:main}
\boxed{\frac{\operatorname{Vol}(\mathcal{P}_n)}{\operatorname{Vol}(\Delta^n)} = \frac{1}{\binom{2n+1}{n}}}.
\end{equation}

Thus, the volume of the inner polytope is exactly $1/\binom{2n+1}{n}$ times the volume of the original simplex. For $n = 1$ (a segment) we get $1/3$, for $n = 2$ — $1/10$, for $n = 3$ — $1/35$ — these are simple cases. For $n = 4$, the formula yields $1/126$, which agrees with Python calculations for the volume of the inner convex hull (volume calculations were performed using both rigorous geometric approaches and the Monte Carlo probabilistic method).

\begin{figure}[htbp]
\centering
Figure~\ref{fig:r2_v4} shows the inner polyhedron for a three-dimensional tetrahedron with trisections edges 1:3:1.

\includegraphics[width=0.45\textwidth]{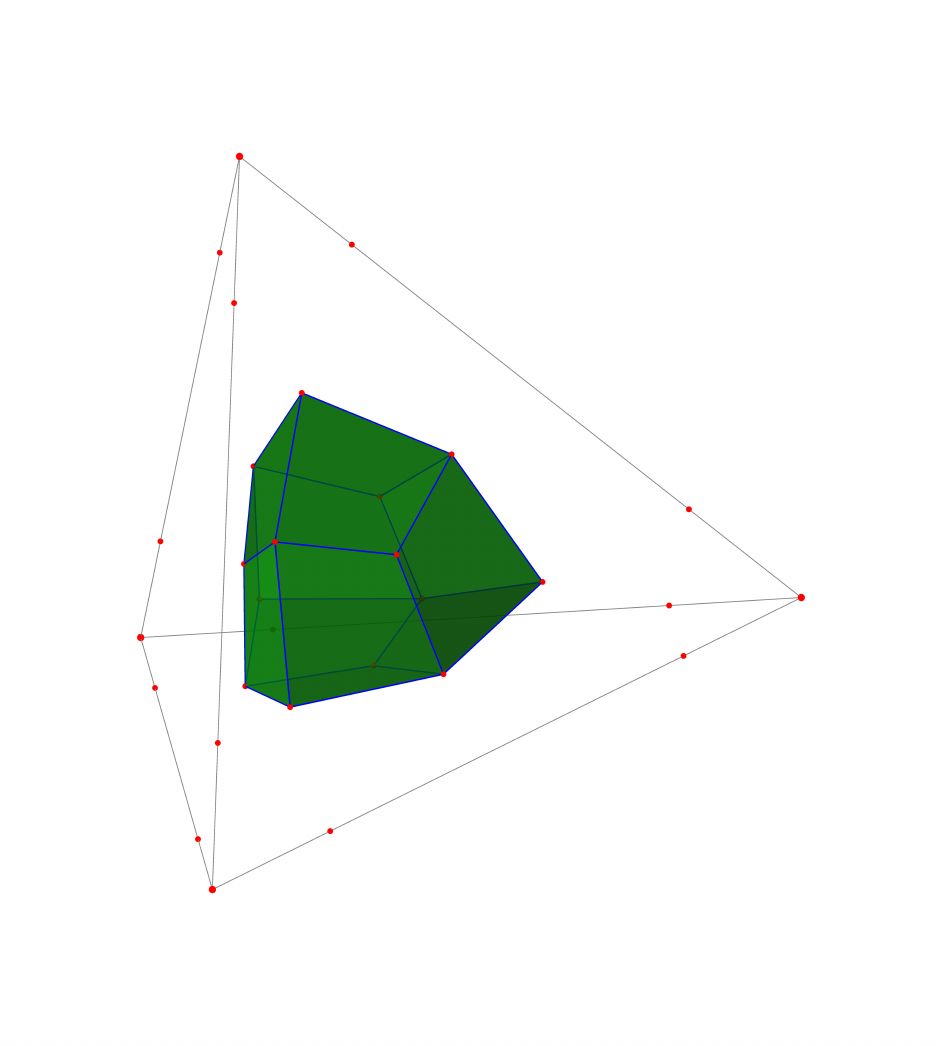}
\caption{Central polytope (deltoidal dodecahedron) obtained in a three-dimensional tetrahedron.}
\label{fig:r2_v4}
\end{figure}

\section{Combinatorial Structure}
In addition to the volume, the inner polytope $\mathcal{P}_n$ has a rich combinatorial structure that can be described through the system of hyperplanes $x_i = 2x_j$.

\subsection{Vertices}
A vertex of $\mathcal{P}_n$ occurs when the maximum number of constraints $x_i \le 2x_j$ is active.
It is easy to verify that the vertices are precisely those points whose barycentric coordinates take only two distinct values, say $\alpha$ and $2\alpha$, with the number of coordinates equal to $\alpha$ being $k$, and the number of coordinates equal to $2\alpha$ being $n+1-k$ ($1\le k \le n$).
From the condition $\sum x_i = 1$, we find $\alpha = 1/(2n+2-k)$ or $\alpha = 1/(n+1+k)$, depending on which value appears $k$ times.
Since the roles of $\alpha$ and $2\alpha$ can be interchanged, each nontrivial partition of the $n+1$ coordinates into two nonempty sets yields two vertices.
Thus, the total number of vertices is
\[
f_0 = 2(2^{n}-1) = 2^{n+1} - 2.
\]

\subsection{Edges}
Two vertices are connected by an edge if their partitions differ by moving exactly one coordinate from the set with value $\alpha$ to the set with value $2\alpha$ (or vice versa) while preserving the overall ratio of values $1/2$ or $2$.
A vertex corresponding to a partition of sizes $(k, n+1-k)$ has degree $k + (n+1-k) = n+1$ for $2 \le k \le n-1$; for $k=1$ or $k=n$, the degree is $n$, since moving a coordinate that would make one of the sets empty is forbidden.
Summing the degrees over all vertices and dividing by $2$, we find the number of edges:
\begin{align*}
f_1 &= \frac{1}{2}\Bigl( 2\cdot n \cdot \binom{n+1}{1} \quad\text{(for $k=1$ and $k=n$)} 
 + \sum_{k=2}^{n-1} 2 \binom{n+1}{k} (n+1) \Bigr) \\
&= (n+1)(2^n - 2).
\end{align*}
The polytope $\mathcal{P}_n$ is a hypersimplex truncated by additional half-spaces $x_i \le 2x_j$, and its graph is a subgraph of the Johnson graph.

\subsection{Faces}
Each hyperface (face of codimension 1) corresponds to an equality $x_i = 2x_j$ for some ordered pair $(i,j)$, $i\neq j$.
There are $n(n+1)$ such ordered pairs, and each inequality $x_i \le 2x_j$ is essential (removing it changes the polytope).
Consequently, the number of hyperfaces is
\[
f_{n-1} = n(n+1).
\]
Each hyperface is combinatorially equivalent to an $(n-1)$-dimensional cube.

The complete $f$-vector can be obtained from the theory of hyperplane arrangements (see~\cite{Zaslavsky}) or via generating functions for face numbers; the polytope $\mathcal{P}_n$ is combinatorially isomorphic to the so-called \emph{shard polytope} associated with $\mathbf{c}$-vectors and dual to the associahedron.

\section{Conclusion and Perspectives}
We have shown that the volume fraction of the central polytope obtained by trisecting all edges of an $n$-dimensional simplex and cutting off the corresponding hyperplanes is exactly
\[
\frac{\operatorname{Vol}(\text{inner polytope})}{\operatorname{Vol}(\text{simplex})} = \frac{1}{\binom{2n+1}{n}}.
\]
This generalizes Marion's theorem for triangles ($n=2$) to tetrahedra ($n=3$) and simplices of higher dimensions.
The proof, based on a probabilistic representation of the uniform distribution on the simplex, transforms the geometric problem into a simple integral and binomial sum.
The same method immediately yields combinatorial data: the number of vertices is $2^{n+1}-2$, edges are $(n+1)(2^n-2)$, and hyperfaces are $n(n+1)$.

Several open questions remain.
One could investigate the volume for an arbitrary integer division ratio $r$ instead of $3$; the probability then leads to a sum of the form $(n+1)\sum_{k=0}^n \binom{n}{k}(-1)^k/(n+1+k(r-1))$, which simplifies to a rational number but, apparently, does not reduce to a single binomial coefficient.
It would also be interesting to study the dual polytope and its volume, or to consider the limiting shape as $n\to\infty$ after appropriate scaling.

\clearpage
\appendix
\section{Probabilistic Method Using the Triangle as an Example}

In this appendix, we reproduce step by step the probabilistic calculation of the area of Marion's central hexagon (the case $n=2$) to illustrate the general method presented in the paper.

\subsection*{Triangle in Barycentric Coordinates}
Consider the standard $2$-simplex (triangle) with vertices $\mathbf{e}_0=(1,0,0)$, $\mathbf{e}_1=(0,1,0)$, $\mathbf{e}_2=(0,0,1)$. Any point of the triangle is given by coordinates $(x_0,x_1,x_2)$ with $x_i\ge 0$, $x_0+x_1+x_2=1$.

Marion's construction: on each side we choose two points dividing it into three equal parts. Through each vertex and the two division points on the opposite side, we draw two lines. All six lines bound the central hexagon $\mathcal{P}_2$. In barycentric coordinates, these lines correspond to the equations
\[
x_1 = 2x_2,\; x_2 = 2x_1 \quad \text{(for edge $\mathbf{e}_0\mathbf{e}_1$ and vertex $\mathbf{e}_2$)},
\]
and similarly for the other edges. The interior region is defined by the inequalities $x_i \le 2x_j$ for all $i\neq j$.

\subsection*{Uniform Distribution on the Triangle}
The key idea: if $Y_0, Y_1, Y_2$ are independent random variables with the standard exponential distribution (density $e^{-y}$, $y\ge 0$), then the vector
\[
(X_0, X_1, X_2) = \left(\frac{Y_0}{S}, \frac{Y_1}{S}, \frac{Y_2}{S}\right), \quad S=Y_0+Y_1+Y_2,
\]
is uniformly distributed in the standard triangle. This is a well-known fact from probability theory (the inverse transform method for the Dirichlet distribution). Therefore, the probability that a random point falls in the hexagon $\mathcal{P}_2$ equals the ratio of areas:
\[
\frac{\operatorname{Area}(\mathcal{P}_2)}{\operatorname{Area}(\Delta^2)} = \mathbb{P}(X_i \le 2X_j \;\forall i\neq j) = \mathbb{P}(Y_i \le 2Y_j \;\forall i,j).
\]

\subsection*{Probability Calculation}
The event $Y_i \le 2Y_j$ for all $i,j$ means that the ratio of the maximum to the minimum of the variables does not exceed $2$. We use symmetry: the probability equals three times the probability that $Y_0$ is the smallest of the three and that $Y_1 \le 2Y_0$, $Y_2 \le 2Y_0$.

Thus, we compute
\[
P = 3 \int_0^\infty f_{Y_0}(y_0) \,\mathbb{P}(Y_1 \le 2y_0 \mid Y_1 \ge y_0) \,\mathbb{P}(Y_2 \le 2y_0 \mid Y_2 \ge y_0) \,dy_0,
\]
where the condition $Y_i \ge y_0$ follows from $y_0$ being the minimum. The density of $Y_0$ is $e^{-y_0}$, and
\[
\mathbb{P}(Y_1 \in [y_0, 2y_0]) = \int_{y_0}^{2y_0} e^{-y} dy = e^{-y_0} - e^{-2y_0}.
\]
Hence,
\[
P = 3 \int_0^\infty e^{-y_0} \bigl(e^{-y_0} - e^{-2y_0}\bigr)^2 \,dy_0.
\]

Transform the integrand:
\[
e^{-y_0} \bigl(e^{-y_0} - e^{-2y_0}\bigr)^2 = e^{-y_0} \cdot e^{-2y_0} (1 - e^{-y_0})^2 = e^{-3y_0} (1 - e^{-y_0})^2.
\]
Thus,
\[
P = 3 \int_0^\infty e^{-3y_0} (1 - e^{-y_0})^2 \,dy_0.
\]

Expand the square: $(1 - e^{-y_0})^2 = 1 - 2e^{-y_0} + e^{-2y_0}$. The integral becomes a sum:
\begin{align*}
P &= 3 \int_0^\infty e^{-3y_0} \bigl(1 - 2e^{-y_0} + e^{-2y_0}\bigr) \,dy_0 \\
&= 3 \left[ \int_0^\infty e^{-3y_0} dy_0 - 2\int_0^\infty e^{-4y_0} dy_0 + \int_0^\infty e^{-5y_0} dy_0 \right] \\
&= 3 \left[ \frac{1}{3} - 2\cdot\frac{1}{4} + \frac{1}{5} \right] \\
&= 3 \left( \frac{1}{3} - \frac{1}{2} + \frac{1}{5} \right) \\
&= 3 \cdot \frac{10 - 15 + 6}{30} = 3 \cdot \frac{1}{30} = \frac{1}{10}.
\end{align*}
Thus, the area of the hexagon is $1/10$ of the area of the original triangle.

\subsection*{Connection with the General Formula}
For $n=2$, the general formula gives $1/\binom{5}{2}=1/10$, which coincides with the obtained result. Note that the sum $1/3 - 2/4 + 1/5$ is equivalent to the expression
$(n+1)\sum_{k=0}^n \frac{\binom{n}{k}(-1)^k}{n+1+k}$ for $n=2$, and after simplification yields the same value.

This example clearly demonstrates how the probabilistic representation reduces a geometric problem to simple integration and algebraic combinatorics. Calculations for higher dimensions proceed analogously, where the power $(1-e^{-y_0})^n$ generates the binomial coefficient $\binom{2n+1}{n}$ in the denominator.

\clearpage
\section{Volume Calculation of the Inner Polytope for a Tetrahedron. Parameters of the Deltoid}

In this appendix, we obtain the central deltoidal dodecahedron for the tetrahedron with vertices $(0,0,0)$, $(1,1,0)$, $(1,0,1)$, $(0,1,1)$ by direct geometric computation, and find its volume and face parameters — congruent deltoids (kites).

\subsection*{Coordinates and Cutting Planes}
Consider the tetrahedron with vertices (in three-dimensional space):
\begin{align*}
A &= (0,0,0), \\
B &= (1,1,0), \\
C &= (1,0,1), \\
D &= (0,1,1).
\end{align*}

\subsection*{Conversion Formula from Barycentric to Cartesian Coordinates}
For this tetrahedron, the transformation matrix from barycentric coordinates $(x_A, x_B, x_C, x_D)$ to Cartesian coordinates $(x, y, z)$ is:
\[
\begin{pmatrix} x \\ y \\ z \end{pmatrix} = 
\begin{pmatrix}
0 & 1 & 1 & 0 \\
0 & 1 & 0 & 1 \\
0 & 0 & 1 & 1
\end{pmatrix}
\cdot
\begin{pmatrix} x_A \\ x_B \\ x_C \\ x_D \end{pmatrix}.
\]

In expanded form:
\[
\begin{cases}
x = x_B + x_C, \\
y = x_B + x_D, \\
z = x_C + x_D.
\end{cases}
\]

\subsection*{Vertices of the Deltoidal Dodecahedron}
According to the general theory, the vertices of the polytope are points where the coordinates take only two values: $\alpha$ and $2\alpha$, with $k\alpha + (4-k)2\alpha = 1$, where $k$ is the number of coordinates equal to $\alpha$. We obtain:
\begin{itemize}
\item $k=1$: $\alpha = 1/7$, $2\alpha = 2/7$ (4 vertices);
\item $k=2$: $\alpha = 1/6$, $2\alpha = 1/3$ (6 vertices);
\item $k=3$: $\alpha = 1/5$, $2\alpha = 2/5$ (4 vertices).
\end{itemize}
Total: $2^4-2 = 14$ vertices.

\subsection*{Table of Vertices of the Deltoidal Dodecahedron}

\begin{table}[htbp]
\centering
\large
\renewcommand{\arraystretch}{1.3}
\begin{tabular}{|c|c|c|c|}
\hline
No. & Letter & Barycentric coordinates & Cartesian coordinates \\
\hline
1 & A & $(2/7, 1/7, 2/7, 2/7)$ & $(3/7, 3/7, 4/7)$ \\
2 & B & $(2/7, 2/7, 1/7, 2/7)$ & $(3/7, 4/7, 3/7)$ \\
3 & C & $(2/7, 2/7, 2/7, 1/7)$ & $(4/7, 3/7, 3/7)$ \\
4 & D & $(1/7, 2/7, 2/7, 2/7)$ & $(4/7, 4/7, 4/7)$ \\
\hline
5 & E & $(1/3, 1/6, 1/6, 1/3)$ & $(1/3, 1/2, 1/2)$ \\
6 & F & $(1/3, 1/6, 1/3, 1/6)$ & $(1/2, 1/3, 1/2)$ \\
7 & G & $(1/3, 1/3, 1/6, 1/6)$ & $(1/2, 1/2, 1/3)$ \\
8 & H & $(1/6, 1/3, 1/6, 1/3)$ & $(1/2, 2/3, 1/2)$ \\
9 & I & $(1/6, 1/3, 1/3, 1/6)$ & $(2/3, 1/2, 1/2)$ \\
10 & J & $(1/6, 1/6, 1/3, 1/3)$ & $(1/2, 1/2, 2/3)$ \\
\hline
11 & K & $(2/5, 1/5, 1/5, 1/5)$ & $(2/5, 2/5, 2/5)$ \\
12 & L & $(1/5, 2/5, 1/5, 1/5)$ & $(3/5, 3/5, 2/5)$ \\
13 & M & $(1/5, 1/5, 2/5, 1/5)$ & $(3/5, 2/5, 3/5)$ \\
14 & N & $(1/5, 1/5, 1/5, 2/5)$ & $(2/5, 3/5, 3/5)$ \\
\hline
\end{tabular}
\caption{Vertices of the deltoidal dodecahedron in barycentric and Cartesian coordinates. 
Vertices 1--4: $k=1$ (4 vertices), 5--10: $k=2$ (6 vertices), 11--14: $k=3$ (4 vertices).}
\label{tab:vertices}
\end{table}

{\small
\textbf{Explanation of the table:}
\begin{itemize}
\item Vertices with $k=1$ (No. 1–4): one coordinate equals $1/7$, three coordinates equal $2/7$.
\item Vertices with $k=2$ (No. 5–10): two coordinates equal $1/6$, two coordinates equal $1/3$.
\item Vertices with $k=3$ (No. 11–14): three coordinates equal $1/5$, one coordinate equals $2/5$.
\item Cartesian coordinates are obtained by $x = x_B + x_C$, $y = x_B + x_D$, $z = x_C + x_D$.
\end{itemize}
}

\subsection*{Volume of the Deltoidal Dodecahedron}
The volume of the original tetrahedron with vertices $(0,0,0)$, $(1,1,0)$, $(1,0,1)$, $(0,1,1)$ is:
\[
V_{\text{tetra}} = \frac{1}{3}.
\]

The volume of the deltoidal dodecahedron can be found by the general formula:
\[
V_{\text{del}} = \frac{V_{\text{tetra}}}{\binom{7}{3}} = \frac{1}{35} \cdot \frac{1}{3} = \frac{1}{105}.
\]

\subsection*{Alternative Calculation of the Volume Fraction via Determinant}
The volume fraction of the inner polytope can also be calculated using the fact that the vertices of the polytope are convex combinations of the simplex vertices. For more details on barycentric coordinates, see the article on MathWorld~\cite{MathWorldBarycentric}.

For the central polytope in a tetrahedron, the volume can be computed via a determinant composed of the barycentric coordinates of four points: three vertices of one face and the center. The determinant
\[
24 \cdot \det \begin{pmatrix}
2/7 & 1/7 & 2/7 & 2/7 \\
1/3 & 1/6 & 1/6 & 1/3 \\
2/5 & 1/5 & 1/5 & 1/5 \\
1/4 & 1/4 & 1/4 & 1/4
\end{pmatrix} = \frac{1}{35}
\]
gives the volume fraction of the inner polytope. The factor $24$ arises because each of the 12 deltoids is divided into two identical triangles, and thus the inner polytope is composed of 24 identical pyramids.

\subsection*{Geometry of the Deltoid}
The faces of the deltoidal dodecahedron are 12 congruent deltoids. Consider one face corresponding to the equality $x_A = 2x_B$. Its vertices (in barycentric coordinates) are:
\begin{align*}
V_1 &= (2/7,\; 1/7,\; 2/7,\; 2/7), \\
V_2 &= (1/3,\; 1/6,\; 1/6,\; 1/3), \\
V_3 &= (1/3,\; 1/6,\; 1/3,\; 1/6), \\
V_4 &= (2/5,\; 1/5,\; 1/5,\; 1/5).
\end{align*}

After converting to Cartesian coordinates using the formulas $x = x_B + x_C$, $y = x_B + x_D$, $z = x_C + x_D$, we obtain:
\begin{align*}
V_1 &\to \left(\frac{3}{7},\;\frac{3}{7},\;\frac{4}{7}\right), \\
V_2 &\to \left(\frac{1}{3},\;\frac{1}{2},\;\frac{1}{2}\right), \\
V_3 &\to \left(\frac{1}{2},\;\frac{1}{3},\;\frac{1}{2}\right), \\
V_4 &\to \left(\frac{2}{5},\;\frac{2}{5},\;\frac{2}{5}\right).
\end{align*}

Computing the side lengths and diagonals:
\begin{itemize}
\item Short sides:
\[
|V_1V_2| = |V_1V_3| = \frac{\sqrt{34}}{42}
\]
\item Long sides:
\[
|V_2V_4| = |V_3V_4| = \frac{\sqrt{22}}{30}
\]
\item Diagonals:
\[
|V_2V_3| = \frac{\sqrt{2}}{6} \text{ (long diagonal)},\qquad 
|V_1V_4| = \frac{\sqrt{38}}{35} \text{ (short diagonal)}.
\]
\end{itemize}

The area of a deltoid is half the product of its diagonals:
\[
S_{\text{del}} = \frac12 \cdot |V_2V_3| \cdot |V_1V_4| 
= \frac12 \cdot \frac{\sqrt{2}}{6} \cdot \frac{\sqrt{38}}{35}
= \frac{\sqrt{76}}{420} = \frac{2\sqrt{19}}{420} = \frac{\sqrt{19}}{210}.
\]

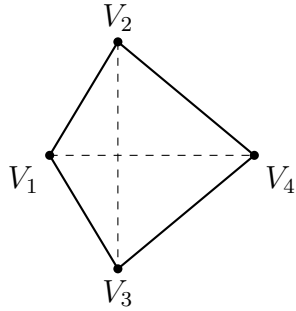
\begin{figure}[htbp]
\centering
\begin{tikzpicture}[scale=6]
\coordinate (V1) at (0,0);
\coordinate (V4) at (0.45,0);
\coordinate (V2) at (0.15,0.25);
\coordinate (V3) at (0.15,-0.25);

\draw[thick] (V1) -- (V2) -- (V4) -- (V3) -- cycle;
\draw[dashed] (V1) -- (V4);
\draw[dashed] (V2) -- (V3);

\fill (V1) circle (0.01) node[below left] {$V_1$};
\fill (V2) circle (0.01) node[above] {$V_2$};
\fill (V3) circle (0.01) node[below] {$V_3$};
\fill (V4) circle (0.01) node[below right] {$V_4$};
\end{tikzpicture}
\caption{Deltoid — a face of the central deltoidal dodecahedron.}
\label{fig:deltoid}
\end{figure}

\subsection*{Volume Verification via Pyramids}
The center of the polytope has barycentric coordinates $(1/4,1/4,1/4,1/4)$, which in Cartesian coordinates gives
\[
O\!\left(\frac12,\frac12,\frac12\right).
\]
The plane of the face $x_A = 2x_B$ in Cartesian coordinates is $3x+3y-z=2$ (derived from the barycentric relations). The distance from the center to the plane:
\[
h = \frac{|3\cdot\frac12+3\cdot\frac12-\frac12-2|}{\sqrt{3^2+3^2+(-1)^2}}
   = \frac{|\frac12|}{\sqrt{19}} = \frac{1}{2\sqrt{19}}.
\]

The volume of the deltoidal dodecahedron, consisting of 12 pyramids with the deltoid base and the vertex at the center, is:
\[
V = 12 \cdot \frac13 S_{\text{del}}\, h = 4 \cdot \frac{\sqrt{19}}{210} \cdot \frac{1}{2\sqrt{19}}
   = \frac{4}{420} = \frac{1}{105}.
\]
This matches the general formula $\displaystyle \frac{V_{\text{tetra}}}{35} = \frac{1/3}{35} = \frac{1}{105}$, confirming the correctness of the calculations.

The code for the deltoidal dodecahedron is given in Appendix C.

\clearpage
\section{Asymptote Code. Central Deltoidal Dodecahedron Inside a Regular Tetrahedron}

Three-dimensional images of such polytopes are conveniently created using the open-source language Asymptote~\cite{Asymptote}. This approach produces high-quality vector interactive images of three-dimensional objects with custom projection and lighting. Below is the code that visualizes the central polytope (deltoidal dodecahedron) obtained in a three-dimensional tetrahedron by trisecting the edges. The reader can copy this code and run it in the Asymptote environment. More information about the language can be found on the official website: \url{https://asymptote.ualberta.ca/}.

\begin{verbatim}
//
// Trisecting edges of tetrahedron a:b:a.
// We can change t = (a+b)/a.
//
real t = 2;

settings.render = 8;
settings.tex = "none";
import three;
size(500);

// Tetrahedron vertices
triple A = (0,0,0);
triple B = (1,1,0);
triple C = (1,0,1);
triple D = (0,1,1);

// Function bit generation
int getBit(int num, int pos) {
    int p = 1;
    for (int k = 0; k < pos; ++k) p *= 2;
    return quotient(num, p) % 2;
}

// 14 Vertex generation
triple[] getVertices(real t) {
    triple[] V;
    for (int mask = 1; mask < 15; ++mask) {
        int k = 0;
        for (int i = 0; i < 4; ++i) {
            if (getBit(mask, i) == 1) ++k;
        }
        if (k == 0 || k == 4) continue;
        real alpha = 1.0 / (k + t * (4 - k));
        real[] bary = new real[4];
        for (int i = 0; i < 4; ++i) {
            if (getBit(mask, i) == 1) bary[i] = alpha;
            else bary[i] = t * alpha;
        }
        V.push(bary[0]*A + bary[1]*B + bary[2]*C + bary[3]*D);
    }
    return V;
}

triple[] V = getVertices(t);

// Edges
int[][] edges = {
    {1,5}, {1,2}, {9,13}, {1,9},
    {2,6}, {5,6}, {2,10},
    {3,5}, {13,11}, {7,11},
    {4,12}, {4,3},
    {5,3}, {5,13},
    {7,9}, {11,12},
    {3,11}, {6,4},
    {8,7}, {8,12},
    {0,4}, {0,2}, {0,8}, {9,10}, {8,10}
};

// Faces construction
int[][] adj;
for (int i = 0; i < 14; ++i) {
    adj[i] = new int[];
}
for (int[] e : edges) {
    adj[e[0]].push(e[1]);
    adj[e[1]].push(e[0]);
}

int[][] faces;
for (int v = 0; v < 14; ++v) {
    for (int u : adj[v]) {
        for (int w : adj[v]) {
            if (u >= w) continue;
            for (int x : adj[u]) {
                if (x == v) continue;
                bool connected = false;
                for (int y : adj[x]) {
                    if (y == w) { connected = true; break; }
                }
                if (connected) {
                    int[] face = {v, u, x, w};
                    int[] sorted = copy(face);
                    int minpos = 0;
                    for (int i = 1; i < 4; ++i) {
                        if (sorted[i] < sorted[minpos]) minpos = i;
                    }
                    int[] rotated = {sorted[minpos], sorted[(minpos+1)%4], sorted[(minpos+2)%4], sorted[(minpos+3)%4]};
                    int[] reversed = {rotated[0], rotated[3], rotated[2], rotated[1]};
                    if (reversed[1] < rotated[1]) rotated = reversed;
                    faces.push(rotated);
                }
            }
        }
    }
}
int[][] uniqueFaces;
for (int[] f : faces) {
    bool exists = false;
    for (int[] uf : uniqueFaces) {
        if (uf[0] == f[0] && uf[1] == f[1] && uf[2] == f[2] && uf[3] == f[3]) {
            exists = true; break;
        }
    }
    if (!exists) uniqueFaces.push(f);
}

// --- Drawing ---
// Tetrahedron edges (gray dashed)
draw(A--B, gray+0.6pt+dashed);
draw(A--C, gray+0.6pt+dashed);
draw(A--D, gray+0.6pt+dashed);
draw(B--C, gray+0.6pt+dashed);
draw(B--D, gray+0.6pt+dashed);
draw(C--D, gray+0.6pt+dashed);

// Polyhedron edges
for (int[] e : edges) {
    draw(V[e[0]] -- V[e[1]], linewidth(0.8pt) + blue);
}

// Face filling
for (int[] f : uniqueFaces) {
    path3 poly = V[f[0]] -- V[f[1]] -- V[f[2]] -- V[f[3]] -- cycle;
    draw(surface(poly), green + opacity(0.7));
    draw(poly, blue + linewidth(1pt));
}

// Polyhedron vertices
for (int i = 0; i < V.length; ++i) {
    dot(V[i], linewidth(4pt) + red);
}

// Tetrahedron vertices (large black)
dot(A, linewidth(4pt) + black);
dot(B, linewidth(4pt) + black);
dot(C, linewidth(4pt) + black);
dot(D, linewidth(4pt) + black);
//label("A", A, SW);
//label("B", B, SE);
//label("C", C, NE);
//label("D", D, NW);

// Division points on edges (small red, no label)
real u1 = 1/(1+t);
real u2 = t/(1+t);
triple[][] tetEdges = {{A,B}, {A,C}, {A,D}, {B,C}, {B,D}, {C,D}};
for (int i = 0; i < tetEdges.length; ++i) {
    triple P = tetEdges[i][0];
    triple Q = tetEdges[i][1];
    dot(P + u1*(Q-P), red+linewidth(4pt));
    dot(P + u2*(Q-P), red+linewidth(4pt));
}

// Optional vertex labels (commented)
//string[] labels = {"V0","V1","V2","V3","V4","V5","V6","V7","V8","V9","V10","V11","V12","V13"};
//for (int i = 0; i < V.length; ++i) {
    //label(labels[i], V[i], align=unit(V[i] - (0.5,0.2,0.1)), fontsize(7pt));
//}

currentprojection = orthographic(2, 1.5, 1);
currentlight = White;
\end{verbatim}

\clearpage
\end{document}